\documentclass[11 pt, a4paper]
{article}

\usepackage[latin1]{inputenc}

\usepackage{amsfonts}
\usepackage{latexsym}
\newcommand{\dis}{\displaystyle}

\newtheorem{Th}[subsection]{Teorema }
\newtheorem{Df}[subsection]{Definici\'on }
\newtheorem{Pro}[subsection]{Proposici\'on }
\newtheorem{remark}[subsection]{Observaci\'on}

\newenvironment{Proof}{\noindent {\bf Demostraci\'on}:}{\QED \smallskip}

\newtheorem{Quest}[subsection]{Cuesti\'on abierta}

\def\topo{topolog\'ia }
\def\topos{topolog\'ias }
\def\eps{\epsilon}
\def\Z{\Bbb Z}
\def\R{\Bbb R}

\def\N{\Bbb N}
\def\T{\Bbb T}

\def\C{\Bbb C}

\def\cont{\mathfrak c}
\newcommand\QED{\hfill QED \medskip}
\begin{document}
\author{D. de la Barrera y E. Mart\'in Peinador }
\title{Estructuras de grupo topol\'ogico en  $(\Z, +)$}




\centerline{\LARGE Estructuras de grupo topol\'ogico en  $(\Z, +)$\footnote{Parcialmente financiado por el MICINN. Proyecto:  MTM2009-14409-C02-01.}}
\vspace{0.5cm}
\centerline{\fontsize{12pt}{1em}\textbf{D. de la Barrera y E. Mart\'in Peinador}}

\centerline{Departamento de Geometr\'{\i}a y Topolog\'{\i}a. Universidad Complutense de Madrid.}

\centerline{\{dbarrera, peinador\}@mat.ucm.es}

\vspace{0.5cm}

\begin{abstract}

{\it En este trabajo hemos considerado la complejidad de las posibles estructuras de grupo topol\'ogico en el grupo de los enteros. Compendiamos algunos resultados nuestros y de otros matem\'aticos para presentar en $\Z$:
\begin{enumerate}

\item Una familia de $2^\cont$ topolog\'ias de grupo precompactas (no compactas) y Hausdorff .

\item Una amplia familia de topolog\'ias de grupo metrizables (no completas, ni precompactas).

\item Una amplia familia de topolog\'ias completas no metrizables.

\end{enumerate}

Dejamos algunas cuestiones abiertas - por ejemplo el problema de la cardinalidad de las familias especificadas en $(2)$ y $(3)$-, que estamos estudiando y  forman parte de la Tesis Doctoral del primer autor.}
\end{abstract}

\centerline{En homenaje al Profesor J.J.  Etayo con afecto y admiraci\'on.}

\section*{\fontsize{12pt}{1em} Introducci\'on}
 \noindent La proverbial frase de Kronecker "los enteros los hizo Dios, los dem\'as n\'umeros son cosa del hombre"\footnote{Die ganzen Zahlen hat der liebe Gott gemacht, alles andere ist Menschenwerk, Kronecker (1823-1891)} quiere resaltar la belleza y la sencillez  de los n\'umeros enteros. Sin embargo, dotar a  su conjunto $\Bbb Z$ de una estructura de  grupo topol\'ogico es una operaci\'on de enorme complejidad. Describimos  a lo largo de este art\'iculo diversos modos de topologizar   $\Bbb Z$ de forma  que   las operaciones suma ordinaria  y  "tomar inverso" sean aplicaciones  continuas, dando lugar as\'i   a  grupos topol\'ogicos con conjunto soporte los n\'umeros enteros.

 Los grupos topol\'ogicos fueron introducidos por Schreier en 1926, en un contexto en el que ya se conoc\'ian los grupos de Lie (grupos continuos de transformaciones).
Hacia los a\~{n}os 30 del siglo pasado hubo una actividad muy fruct\'ifera en grupos topol\'ogicos: de esa \'epoca datan los trabajos de Weyl, la introducci\'on de la medida de Haar y los trabajos de Pontryagin estableciendo el Teorema de la dualidad de grupos topol\'ogicos localmente compactos abelianos, que es la piedra angular para el an\'alisis arm\'onico. Un grupo importante en este contexto  es el toro $\T$, grupo de los complejos m\'odulo uno, con la topolog\'ia inducida por la usual de $\C$. En efecto, $\T$ es el objeto dualizante en la teor\'ia de dualidad de Pontryagin. Los homomorfismos de un grupo abeliano  cualquiera $G$ en $\T$ se denominan caracteres, y su conjunto $Hom (G, \T)$ constituye un grupo respecto de la operaci\'on definida puntualmente.

Damos en primer lugar la definici\'on formal y propiedades elementales de los grupos topol\'ogicos que nos permitan c\'omodamente hilvanar nuestros argumentos sobre $\Bbb Z$.
\begin{Df}\label{1}
Un grupo topol\'ogico es una terna $(G, *, \tau)$
d\'onde $G$ designa un conjunto, $*:G\times G \to G$ una operaci\'on binaria que da a $G$ estructura de grupo y $\tau$ una topolog\'ia en $G$  que hace continua la operaci\'on $*$ -considerando en $G \times G$ la topolog\'ia producto-, as\'i como la inversi\'on $G \to G$ definida por $ x \in G \mapsto x^{-1}$.
\end{Df}

De la propia definici\'on surge el corolario de que las traslaciones a derecha e izquierda:
$$ r_a, l_a: G \to G
 \mbox{  definidas por  } r_a (x)=  x*a \mbox{  y  } l_a(x) = a* x $$
 son homeomorfismos, para cada $a \in G$. En consecuencia los entornos de un punto cualquiern  $a \in G$ est\'an perfectamente determinados por los entornos del elemento neutro $e$, mediante la correspondiente traslaci\'on.

La \topo $\tau$ de la definici\'on \ref{1} se dir\'a que es una \topo de grupo para $G$, o tambi\'en que la \topo $\tau$ es compatible con la estructura de grupo de $G$. La estructura algebraica de grupo interacciona con la estructura topol\'ogica, y algunas propiedades topol\'ogicas son autom\'aticas, mientras que otras se ven notablemente  reforzadas. Por ejemplo, todo grupo topol\'ogico es un espacio topol\'ogico homog\'eneo,
 completamente regular  (y por tanto uniformizable \cite[(VII.7.22)]{Out2}).  Las  proposiciones siguientes  
   son tambi\'en muestra de \'ello.

\begin{Pro}\label{Hd}
Todo grupo topol\'ogico $T_1$ es tambi\'en $T_2$ (es decir, Hausdorff).
\end{Pro}
\begin{Pro} \label{metr}
Todo grupo topol\'ogico I-numerable y Hausdorff es metrizable (Teorema de Birkhoff-Kakutani).
\end{Pro}

En un conjunto cualquiera $X$ pueden definirse  topolog\'ias sencillas  relacionadas con su \'algebra de Boole $\mathcal{ P}(X)$, c\'omo la topolog\'ia discreta $\tau_d$ \'o la cofinita $\tau_c$. Si en el lugar del conjunto $X$ tomamos un grupo algebraico $(G, *)$ cabe plantearse si dichas topolog\'ias son compatibles con la estructura de grupo. 

\begin{Pro}
Sea $G$ un grupo infinito y sean $\tau_d$ y $\tau_c$ las topolog\'ias   en $G$ discreta  y cofinita respectivamente. Entonces $(G, \tau_d)$ es un grupo topol\'ogico, mientras que $(G, \tau_c)$ no lo es.
\end{Pro}

\begin{Proof}
La primera afirmaci\'on es evidente, mientras que la segunda es consecuencia de la Proposi\'on \ref {Hd}, puesto que $G$ dotado de $\tau_c$ es un espacio $T_1$, pero no es $T_2$.
\end{Proof}

 En lo que sigue denominamos simplemente $G$ a un grupo algebraico o topol\'ogico, la operaci\'on y -en su caso- la topolog\'ia se sobreentienden, y en lugar de $x*y$ escribimos $xy$.

Bas\'andonos en la Proposici\'on \ref{Hd}, s\'olo nos van a interesar los grupos topol\'ogicos de Hausdorff, puesto que una m\'inima propiedad de separaci\'on (incluso $T_0$) en un grupo topol\'ogico ya implica que se cumple el axioma $T_2$. A partir de ahora los grupos considerados en este trabajo son {\bf abelianos y Hausdorff}.

{\bf Sistema de entornos del elemento neutro $e \in G$}.

Describimos las propiedades esenciales  de una subfamilia $\mathcal{N}(e)$ de partes de $G$, para poder asegurar que constituyen un sistema  de entornos del elemento neutro para una topolog\'ia  de grupo en $G$.

\begin{itemize}
\item[(G1)]  $e \in V, \forall V \in \mathcal{N}(e)$
\item[(G2)] $V_1  \cap   V_2 \in \mathcal{N}(e)$, siempre que $V_1, V_2 \in \mathcal{N}(e)$
 \item[(G3)] Si para $W \subset G$ existe  $V \in \mathcal{N}(e)$  tal que $V \subseteq W$,
entonces $ W \in \mathcal{N}(e)$
\item[(G4)] $\forall W \in \mathcal{N}(e)$, existe $V\in \mathcal{N}(e)$ tal que $VV: = \{xy: x, y \in V \} \subset W$
\item[(G5)] $ W \in \mathcal{N}(e) \Rightarrow  W^{-1} : = \{x^{-1}: x \in W \} \in \mathcal{N}(e)$
\end{itemize}

\section{\fontsize{12pt}{1em} Topolog\'ias p-\'adicas}

\noindent Sea $p$ un n\'umero primo. La familia $$\mathcal{B} = \{p^n \Bbb Z: n \in \Bbb N\}$$
es base de entornos de $0$ para una topolog\'ia de grupo en $\Bbb Z$, que se denomina la topolog\'ia p-\'adica,  y que denotaremos por $\lambda_p$. F\'acilmente se comprueba que $\mathcal{B}$ cumple las propiedades $G1, G2, G4, G5,$ y adem\'as   los conjuntos de $\mathcal{B}$ son subgrupos de $\Bbb Z$. En general, una topolog\'ia de grupo que admite una base de entornos de cero formada por subgrupos se denomina una {\it topolog\'ia lineal}.\\
La simple observaci\'on de que el conjunto de los n\'umeros primos $\Bbb P$ tiene cardinal $\aleph_0$  nos permite afirmar que las topolog\'ias p-\'adicas constituyen una  familia infinita de topolog\'ias en $\Bbb Z$. En efecto, si $p, q$ son primos distintos
$(\Bbb Z, \lambda_p)$ y $(\Bbb Z, \lambda_q)$ no son topol\'ogicamente isomorfos. Los \'unicos automorfismos de $\Bbb Z$ son la identidad y la inversi\'on, que {\bf no} transforman entornos de cero en $\lambda_p$ en entornos de cero en $\lambda_q$.

Podemos pensar qu\'e otras \topos lineales existen en $\Z$.   C\'omo  los subgrupos de $\Z$ son de la forma $m\Z$, para un entero $m$,   las \topos lineales en $\Z$ se caracterizan f\'acilmente a trav\'es de las llamadas $D$-sucesiones, que  juegan el mismo papel que
 la sucesi\'on $p, p^2, \dots p^n, \dots$  en la \topo $p$-\'adica.
 \begin{Df}
  Una sucesi\'on $b_1, b_2, \dots$ de n\'umeros naturales se dir\' a que es una $D$-sucesi\'on si:

\begin{itemize}

\item $b_1 = 1$.

\item $b_n\neq b_{n+1}$ para todo n\'umero natural $n$.

\item $b_{n+1}$ es m\'ultiplo de $b_n$ para todo n\'umero natural $n$.

\end{itemize}
\end{Df}

Una $D$-sucesi\'on ${\bf b} = \{b_n, n\in \N \}$ define una  \topo lineal  en $\Z$, que   denominaremos  \topo $b$-\'adica $\lambda _b$, mediante la siguiente base de entornos de cero: $$\mathcal{B} = \{b_n \Bbb Z: n \in \Bbb N\}$$
Todas las \topos lineales no discretas de Hausdorff en $\Z$ se definen mediante una $D$-sucesi\'on (v. \cite[(2,1)]{AusBar2012}) .
Destacamos el siguiente hecho f\'acil de demostrar:
\begin{Pro}
Sea ${\bf b} = \{b_n, n\in \N \}$ una  $D$-sucesi\'on.
En $(\Z, \lambda_b)$  la sucesi\'on 
${\bf b}$ converge a $0$. En particular, $p^n \to 0 $ en la \topo $p$-\'adica
\end{Pro}

Teniendo en cuenta que un subgrupo abierto de un grupo topol\'ogico es tambi\'en cerrado obtenemos la siguiente propiedad:

\begin{Pro}
Un grupo lineal $G$ es cero-dimensional y por tanto es totalmente inconexo. En particular  $(\Z, \lambda_b)$  es cero dimensional, para cualquier  $D$-sucesi\'on ${\bf b}$.

\end{Pro}

\section{\fontsize{12pt}{1em} Topolog\'ias precompactas}

\noindent Las topolog\'ias $b$-\'adicas en $\Z$ descritas en la anterior secci\'on son precompactas, en el siguiente sentido:

\begin{Df}
Se dir\'a que  un grupo topol\'ogico $(G, \tau)$ es  precompacto si para todo $V$ entorno de cero existe un subconjunto finito $F$ tal que $G = F + V$, es decir $G$ se obtiene como uni\'on finita de trasladados de $V$.
\end{Df}

Es una consecuencia directa de la definici\'on anterior que todo grupo compacto es tambi\'en precompacto. Por el teorema de categor\'ia de Baire, un grupo infinito compacto y $T_2$, o metrizable y completo no puede ser {\bf numerable}. Por tanto el grupo $\Z$ no admite ninguna \topo de Hausdorff compacta. Sin embargo veremos que admite $2^{\cont}$ \topos precompactas y $T_2$. Por otro lado, no puede haber una topolog\'ia de grupo metrizable y completa en $\Z$. En la secci\'on $\ref{secciontopologiascompletas}$ daremos una amplia familia de topolog\'ias de grupo completas en $\Z$ (no metrizables).

Las \topos precompactas en un grupo abeliano $G$ est\'an relacionadas con el grupo de los caracteres $Hom (G, \T)$.  Destacamos la definici\'on de grupo dual que usaremos en lo sucesivo.

\begin{Df}
Sea $(G, \tau)$ un grupo topol\'ogico. El grupo dual de $G$ que designaremos por $G^\wedge$ es el grupo de los caracteres continuos, es decir $G^\wedge \leq Hom (G, \T)$. Si $\tau $ es la \topo discreta $G^\wedge = Hom (G, \T)$.
\end{Df}

Si $G$ es el grupo de los enteros $\Z$, cualquier homomorfismo de  $\Z$ en $\T$ viene definido por la imagen de su generador $1$, y su conjunto $Hom (\Z, \T) $,  puede identificarse con $\T$. La correspondencia $\varphi\mapsto\varphi(1)$ es de hecho un isomorfismo entre los grupos $Hom(\Z,\T)$ y $\T$.  M\'as a\'un, es un isomorfismo topol\'ogico si $Hom(\Z,\T)$ est\'a dotado de la topolog\'ia de la convergencia puntal y $\T$ tiene  la \topo eucl\'idea ordinaria. Por tanto el dual de $\Z$ con cualquier topolog\'ia de grupo es un subgrupo de $\T$.

Se demuestra f\'acilmente que el grupo dual de $(\Z, \lambda_p)$ es precisamente el grupo de Pr$\ddot{u}$fer $\Z(\mathbf{p}^\infty)$, formado por todas las raices $p^n$-simas de la unidad de $\T$, para $n \in \N$ cualquiera.

El siguiente resultado, consecuencia del teorema de Peter-Weyl permite entender que los grupos precompactos tienen la \topo d\'ebil correspondiente a la familia de sus caracteres continuos.
\begin{Th}\cite[22.14]{HR63} Si $(G, \tau)$ es un grupo compacto y $T_2$, existen suficientes caracteres continuos para separar puntos de $G$. Esto es, para cada par $g_1, g_2 \in G$ con $g_1 \neq g_2$ existe $\varphi \in G^\wedge$ tal que $\varphi(g_1) \neq \varphi(g_2)$.
\end{Th}

 Este teorema precisamente implica que todo grupo compacto y Hausdorff $G$ se encaja en el producto $\T^{G^\wedge}$ 
  mediante un isomorfismo topol\'ogico (no suprayectivo en general). Adem\'as la \topo original de $G$ es la d\'ebil correspondiente  a su familia de caracteres continuos.

  Es inmediato probar que un subgrupo de un grupo compacto es precompacto, y por otra parte un grupo precompacto se sumerge de forma est\'andar en su "completado"  que es compacto. Por tanto, la clase formada por todos los grupos precompactos se pueden identificar a la clase formada por   los subgrupos de los  grupos compactos. Los siguientes  teoremas, extraidos de \cite[(1.2)]{ComRos1964}, nos proporcionar\'an un m\'etodo para construir todas las \topos precompactas de un grupo $G$ y en particular de $\Z$:

 \begin{Th}\label{CR}
 Sea $(G, \tau)$  un grupo topol\'ogico y $G^\wedge$ su grupo dual. Las siguientes afirmaciones son equivalentes:\\
 a) $(G, \tau)$  es precompacto.\\
 b) La \topo $\tau$ coincide con la \topo d\'ebil en $G$ correspondiente a la familia de homomorfismos $G^\wedge$.
 \end{Th}

 \begin{Th}\label{CET}
 Si $G$ es un grupo abeliano, $H$ un subgrupo de $Hom (G, \T)$ que separa puntos de $G$, y $\tau_H$ la \topo d\'ebil en $G$ correspondiente a la familia $H$, el  grupo dual $(G, \tau_H)^\wedge$ es precisamente $H$.
 \end{Th}

 \begin{Pro}\label{proposicion2.6}

 Si $H$ es un subgrupo infinito de $\T$, entonces $H$ separa puntos de $\Z$

 \end{Pro}

 \begin{Proof}
 En efecto, si $m $ es entero no nulo, existe alg\'un $\alpha \in H$ con $\alpha^m \neq 1$, ya que s\'olo hay $m$ ra\'ices $m$-simas de 1 y $H$ es infinito.

 \end{Proof}

 Ahora estamos en condiciones de construir todas las topolog\'ias precompactas de Hausdorff en $\Z$. Por los Teoremas \ref{CR} y \ref{CET}, se justifica el siguiente procedimiento. Sea $H$ un subgrupo infinito del c\'irculo unidad del plano complejo. 
    Definimos la topolog\'ia $\tau_H$ c\'omo la topolog\'ia d\'ebil asociada a $H$; es decir, la que tiene como subbase de entornos de cero la familia  $\mathcal{B}=\{f^{-1}(V):V\in\mathcal{N}_0(\Bbb T) \mbox{ y }f\in H\}$. Se  comprueba directamente que $\tau_H$ es una topolog\'ia de grupo.

Dependiendo del subgrupo $H$ que se tome, se pueden obtener diferentes propiedades de la topolog\'ia $\tau_H$. Por ejemplo:

\begin{itemize}

\item Si $H_1<H_2$, entonces $\tau_{H_1}<\tau_{H_2}$.

\item El conjunto de subgrupos de $\Bbb T$ y el conjunto de topolog\'ias precompactas en $\Bbb Z$ est\'an en correspondencia biyectiva.

\item El subgrupo $H$ es numerable si y s\'olo si $\tau_H$ es metrizable.

\end{itemize}

Denotamos por $\mathfrak{P}$ al conjunto de topolog\'ias precompactas en $\Z$.

\begin{Pro}\label{muchas}

Existen $2^{\cont}$ topolog\'ias precompactas y Hausdorff en $\Z$.

\end{Pro}
\begin{Proof}
Sea $\mathfrak{B}$ una base de Hamel de $\R$, tal que $1\in\mathfrak{B} $. Designamos por $\alpha_b=e^{2\pi i b}$, con $b\in\mathfrak{B} $. Para cada $b\neq 1$, el grupo $\langle \alpha_b\rangle$ generado por $\alpha_b$ es un subgrupo infinito de $\T$ y por la proposici\'on \ref{proposicion2.6} separa puntos de $\Z$. Para cualquier $M\subseteq\mathfrak{B} $, el grupo $H_M$ engendrado por $\{\alpha_m: m\in M\}$ da lugar a una topolog\'ia precompacta y Hausdorff en $\Z$. Teniendo en cuenta que $card(\mathcal{P(\mathfrak{B} )})=2^\cont$ y que si $M\neq M'$, las topolog\'ias $\tau_{H_M}$ y $\tau_{H_{M'}}$ son distintas, obtenemos que $card({\mathfrak{P}})\geq2^\cont$. Por otra parte en $\Z$ no puede haber m\'as de $2^\cont$ topolog\'ias.

\end{Proof}
\begin{remark} {\em Aunque hemos dado una demostraci\'on  directa de la Proposici\'on  \ref{muchas}, el resultado  ya era conocido. De hecho en \cite{Rac2002} se prueba  que  hay dos familias  ${\mathcal A}_1$ y ${\mathcal A}_2$ de $2^\cont$ topolog\'ias precompactas  en $\Z$ -no homeomorfas dos a dos-, tales que las \topos de ${\mathcal A}_1$  no tienen  sucesiones convergentes no triviales, mientras que las de ${\mathcal A}_2$  si las tienen.  Teniendo en cuenta  que toda topolog\'ia   metrizable en cualquier conjunto  infinito  tiene   sucesiones  convergentes no triviales, las \topos de ${\mathcal A}_1$  no son metrizables. }
\end{remark}

La familia $\mathfrak{P}$ tiene como m\'aximo la topolog\'ia $\tau_\T$ o topolog\'ia d\'ebil inducida por todos los caracteres de $\Z$. Siguiendo a Van Douwen en su magn\'ifico trabajo \cite{VD1990}, un grupo abeliano discreto dotado de la topolog\'ia d\'ebil asociada a $Hom(G,\T)$ se denota por $G^{\#}$. En dicho trabajo, se da una prueba de que para cualquier grupo abeliano infinito $G$, $G^{\#}$ no tiene sucesiones convergentes no triviales. En \cite[lema 1]{BanMar1996}, se da -con otros fines- una prueba directa de este hecho para $\Z^\#$.

El cl\'asico Teorema de Pontryagin-Van Kampen sugiere  que la topolog\'ia natural en un grupo dual debe ser la compacto-abierta.  La noci\'on de reflexividad tambi\'en se cimienta en dicho teorema: se dir\'a que un grupo topol\'ogico $G$ es {\it reflexivo} si la evaluaci\'on can\'onica $\alpha_G$ de $G$ en $G^{\wedge\wedge}$ es isomorfismo topol\'ogico, donde ambos $G^\wedge$ y $G^{\wedge\wedge}$ est\'an dotados de la correspondiente topolog\'ia compacto abierta. La clase de los grupos reflexivos abarca los localmente compactos y Hausdorff, precisamente   \'esta es la afirmaci\'on del famoso Teorema mencionado. Ya en los 50 del siglo pasado se  obtuvieron  otros grupos reflexivos no localmente compactos. Recientemente,
Gabriyelyan ha encontrado una topolog\'ia reflexiva no discreta  en el grupo de los enteros, \cite{Gab2010}. Esto es un hecho asombroso,  sin embargo   pensamos  que  la siguiente   pregunta -por ser a\'un m\'as exigente-  tendr\'a una soluci\'on negativa:

\begin{Quest}
( Tkachenko, 2009) \textquestiondown Existe alguna topolog\'ia precompacta $\nu$ en $\Z$ tal que $(\Z,\nu)$ sea reflexivo?
\end{Quest}

\section{\fontsize{12pt}{1em} Topolog\'ias de convergencia uniforme}

\noindent En \cite{Bou1966} (p\'agina 24, ejercicio 2), Bourbaki afirma que las  topolog\'ias  no discretas en el grupo de los enteros  son precompactas. Sin embargo, esta afirmaci\'on es err\'onea, como veremos al final de esta secci\'on.

Vamos a definir mediante una pseudom\'etrica en $\Z$ una topolog\'ia de convergencia uniforme en un subconjunto $S\subseteq\T$.

\begin{Df}

Fijamos un subconjunto $S\subset \T$. La aplicacion $d_S:\Z^2\rightarrow \R_+$  definida de la siguiente manera: $\dis d_S(m,n)=\sup_{\alpha\in S}|\alpha^m-\alpha^n|$
es una seudom\'etrica en $\Z$ acotada por $2$.

\end{Df}

\begin{Pro}

La pseudom\'etrica $d_S$ define una topolog\'ia de grupo en $\Z$ que denotaremos por $\rho_S$. Si $S$ separa puntos de $\Z$, entonces $d_S$ es una m\'etrica.

\end{Pro}

\begin{Proof}
Toda pseudom\'etrica invariante por traslaciones da lugar a una topolog\'ia de grupo \cite{Kl} (tomando como entornos de cero las bolas $B_n=\{z\in\Z:d_S(z,0)<\frac{1}{n}\}$). Es claro que $d_S$ es invariante por traslaciones ya que $\dis d_S(m+k,n+k)= \sup |\alpha^{m+k}-\alpha^{n+k}|=\sup|\alpha^k||\alpha^m-\alpha^n|= \sup|\alpha^m-\alpha^n|=d_S(m,n)$.

\end{Proof}
\begin{remark}

La topolog\'ia $\rho_S$, es la topolog\'ia de convergencia uniforme en el conjunto $S$.  Con notaci\'on aditiva para $\T$, es decir considerado $\T$ c\'omo el grupo cociente $\R / \Z$, se comprueba en   \cite{AusBar2012} que los conjuntos $V_{S,n}:=\{k\in\Bbb Z :  kx+\Bbb Z \in[-\frac{1}{4n},\frac{1}{4n}]+\Bbb Z \mbox{ para todo } x\in S\}$ con $n\in\N$ constituyen una base de entornos de $0$ para $\rho_S$.
\end{remark}

\begin{remark}

Sean $A\subseteq B\subseteq\T$, entonces $\rho_A\leq\rho_B$.

\end{remark}

\begin{Pro}

Si $S$ es denso, entonces $d_S$ da lugar a la topolog\'ia discreta.

\end{Pro}

\begin{Proof}
Vamos a ver que fijados $m,n\in\Z$, tenemos que $d_S(m,n)=d_\T(m,n)$. Para ello tomamos $\alpha\in\T$. Sea $\eps>0$. Por ser $S$ denso, existe $s\in S$ tal que $|s-\alpha|<\frac{\eps}{|m|+|n|}$. Ahora bien, $s^m-\alpha^m=(s-\alpha)(s^{m-1}+s^{m-2}\alpha+\cdots+\alpha^{m-1})$.

Tomando m\'odulos, tenemos que $|s^m-\alpha^m|=|s-\alpha||s^{m-1}+s^{m-2}\alpha+\cdots+\alpha^{m-1}|<\frac{\eps|m|}{|m|+|n|}$.

Por otro lado $\alpha^m-\alpha^n=\alpha^m-s^m+s^m-s^n+s^n-\alpha^n$. Es decir, $|\alpha^m-\alpha^n|\leq|\alpha^m-s^m|+|s^m-s^n|+|s^n-\alpha^n|<\frac{\eps|m|}{|m|+|n|}+d_S(m,n)+\frac{\eps|n|}{|m|+|n|}=d_S(m,n)+\eps$. Tomando supremos tenemos que $d_S(m,n)\leq d_\T(m,n)< d_S(m,n)+\eps$. Al ser la desigualdad cierta para todo $\epsilon >0$, se deduce  que $d_A(m,n)=d_\T(m,n)$.

\end{Proof}

La densidad de $S$ {\bf no} es necesaria para que la topolog\'ia $\rho_S$ sea discreta. Sin embargo, hay "muchos" subconjuntos $S\subseteq \T$ que dan lugar a topolog\'ias uniformes no discretas. En \cite{AusBar2012}  se da una familia de sucesiones $S\subseteq \T$ que verifican:\\ 1)  $(\Z,\rho_S)^\wedge=\langle S\rangle$ y\\ 2)   la \topo  precompacta asociada a $\langle S\rangle$, $\tau_{\langle S \rangle}$  no coincide con $\rho_S$.

Claramente 1) implica que $\rho_S$ no es discreta. Adem\'as,  el dual de $(\Z, \tau_{\langle S \rangle})$ es   $
\langle S \rangle$ (por \ref{CET}), y as\'i los duales de  $(\Z, \tau_{\langle S \rangle})$  y de  $(\Z,\rho_S)$ coinciden. Teniendo en cuenta el Teorema \ref{CR}, s\'olo puede haber una topolog\'ia precompacta en el grupo  $\Z$ entre todas aqu\'ellas que dan lugar al mismo dual. Por 2) obtenemos que $ \rho_S $ no es precompacta.

Cada $ \rho_S $ definida con el criterio anterior, es   una topolog\'ia de grupo {\bf no precompacta} en $\Z$, lo que  contradice la afirmaci\'on de Bourbaki arriba mencionada.

La siguiente proposici\'on se demuestra directamente:

\begin{Pro}

Para $S\subseteq \T$ se obtiene $S\subseteq (\Z,\rho_S)^\wedge$.

\end{Pro}

\begin{Pro}
 La topolog\'ia $\rho_S$ no es  precompacta para ning\'un subconjunto infinito $S\subseteq\T$.
 \end{Pro}

 \begin{Proof}
 Procedemos por reducci\'on al absurdo. Por el teorema \ref{CR}, tenemos que en el caso de que $\rho_S$ fuera precompacta, en particular, ser\'ia la topolog\'ia d\'ebil asociada al subgrupo $H=(\Z,\rho_S)^\wedge\leq\T$; es decir, $\tau_H$. Es claro que $\tau_H\leq\rho_S$. Ahora bien, supongamos que $\rho_S\leq\tau_H$;  para cada entorno de $0$ $U_\rho$ de $\rho_S$, deber\'iamos tener otro entorno de $0$ $U_\tau$ de $\tau_H$ de tal manera que $U_\tau\subseteq U_\rho$.

 Los entornos de $\tau_H$ son de la forma $U_F:=\{k\in\Z|\varphi(k)\in\T_+\mbox{ para todo }\varphi\in F \mbox{ con } F\subseteq(\Z,\rho_S)^\wedge\mbox{ finito}\}$. Por otra parte, $U_S:=\{k\in\Z|\varphi(k)\in\T_+\mbox{ para todo }\varphi\in S\}$ es un entorno de $\rho_S$. Supongamos que existe $F \subset H$ tal que $U_F\subseteq U_S$. Ahora  consideramos $V_F:=\{\varphi\in(\Z,\rho_S)^\wedge|\varphi(x)\in\T_+\mbox{ para todo }x\in U_F\}$ y an\'alogamente, $V_S:=\{\varphi\in(\Z,\rho_S)^\wedge|\varphi(x)\in\T_+\mbox{ para todo }x\in U_S\}$. Por las propias definiciones y la hip\'otesis de que $\rho_S=\tau_H$, tenemos que $V_S\subseteq V_F$ y adem\'as $S\subseteq V_S$. De acuerdo con \cite[7.11]{tesislydia}   $V_F$ es finito, por lo cu\'al es imposible que $V_S\subseteq V_F$. Consecuentemente, $\rho_S\neq \tau_H$.

\end{Proof}

\begin{Quest}
\begin{enumerate}
\item Computar $(\Z,\rho_S)^\wedge$ para $S\subseteq\T$ cualquiera.

\item Dar condiciones necesarias y suficientes para que $\rho_S\neq \rho_R$ con $S,R\subseteq\T$.

\item Calcular la cardinalidad de la familia $\{\rho_S: S\subseteq\T\}$
\end{enumerate}
\end{Quest}
\section{\fontsize{12pt}{1em} Topolog\'ias completas en $\Z$}\label{secciontopologiascompletas}

\noindent Hasta ahora hemos mostrado ejemplos de topolog\'ias  en $\Z$ -metrizables o no-,  que en ning\'un caso son completas.
En esta secci\'on   presentamos una familia de topolog\'ias  completas no discretas, que    debido al teorema de categor\'ia de Baire no son metrizables (ni localmente compactas).

En \cite{ProZel1999} Protasov y Zelenyuk se cuestionan lo siguiente: dada una sucesi\'on $a=(a_n)$ en un grupo $G$, \textquestiondown existe una topolog\'ia de grupo en G que sea de Hausdorff  y tal  que la sucesion $a$ sea convergente a  cero en dicha topolog\'ia? La respuesta depende de la sucesi\'on. Por ejemplo la sucesion $(a_n)=(n^2)\subseteq\Z$ no converge  a $0$ para ninguna topolog\'ia de grupo y Hausdorff en $\Z$. As\'i surge la siguiente noci\'on.

\begin{Df}

Sea $G$ un grupo y $a=(a_n)\subset G$, una sucesi\'on. Diremos que $a$ es una $T$- sucesi\'on si existe una topolog\'ia de grupo y de Hausdorff $\tau$ en $G$ tal que $a_n\rightarrow 0_G$ en $\tau$, siendo $0_G$ el elemento neutro de $G$.

\end{Df}

 Mediante una sencilla aplicaci\'on del lema de Zorn se prueba que si $a\subset G$ es una $T$-sucesi\'on, entonces  existe una topolog\'ia en $G$ de grupo m\'as fina  entre todas aqu\'ellas que      hacen    nula  la sucesi\'on $a =(a_n)$. Protasov y Zelenyuk usan el t\'ermino {\it topolog\'ia determinada por la sucesi\'on } $a$ para designarla y hacen una construcci\'on de la misma que informalmente consiste en: definir primero unos subconjuntos de $G$ que incluyen las colas de la sucesi\'on $a =(a_n)$, y tomar despu\'es todas las sumas finitas de elementos de los subconjuntos anteriores. El proceso da lugar a una base de entornos del neutro para una topolog\'ia de grupo en $G$ que es la m\'as fina entre todas las que hacen nula la sucesi\'on $a$. Lo escribimos formalmente como sigue:

\begin{Df}[\cite{ProZel1999}]

Sea $G$ un grupo, sea $a=(a_n)$ una $T$-sucesi\'on en $G$ y $(n_i)_{i\in \N}$ una sucesi\'on de naturales. Definimos

\begin{itemize}

\item $A^*_m:=\{\pm a_n|n\geq m\}\cup\{0_G\}$.

\item $A(k,m):=\{g_0+\cdots +g_k| g_i\in A^*_m\, i\in\{0,\dots ,k\}\}$.

\item $[n_1,\dots, n_k]:=\{g_1+\cdots + g_k: g_i\in A^*_{n_i}, i=1,\dots, k\}$.

\item $\dis V_{(n_i)}=\bigcup^\infty_{k=1}[n_1,\dots,n_k]$.

\end{itemize}
\end{Df}

\begin{Pro}

La familia $\{V_{(n_i)}:(n_i)\in \N^\N\}$ es una base de entornos de $0_G$ para una topolog\'ia de grupo $\mathcal{T}_{\{a_n\}}$ en $G$, que es la m\'as fina de todas aqu\'ellas en las que convege $a$. El s\'imbolo  $G_{\{a_n\}}$ designar\'a al grupo $G$ dotado de  $\mathcal{T}_{\{a_n\}}$.

\end{Pro}

La topolog\'ia $\mathcal{T}_{\{a_n\}}$ es necesariamente de Hausdorff, por la definici\'on de $T$-sucesi\'on. El siguiente teorema, cuya demostraci\'on no incluimos por ser muy t\'ecnica, es el resultado m\'as importante de la secci\'on:

\begin{Th}\cite[Theorem 8]{ProZel1999}

Sea $G$ un grupo y sea $a=(a_n)$ una $T$-sucesi\'on. Entonces la topolog\'ia $\mathcal{T}_{\{a_n\}}$ es completa.

\end{Th}

Para una $D$-sucesi\'on $ a$ en $\Z$ es claro que se puede construir la topolog\'ia $\mathcal{T}_{\{a_n\}}$ en $\Z$, que de acuerdo con el teorema anterior, es completa y no metrizable.
Las siguientes propiedades de $\Z_{\{a_n\}}$  tienen especial inter\'es: la primera porque resuelve un problema de Malykhin de varios a\~{n}os de antig\"{u}edad y la segunda  nos permite obtener   que el grupo dual es metrizable y $ k$-espacio.

\begin{remark}\label{kW}
\begin{itemize}

\item[(1)]   $\Z_{\{a_n\}}$ es un grupo secuencial, pero no es Frechet-Urysohn (ver \cite[Theorem 7 y Theorem 6]{ProZel1999}).

\item[(2)] $\Z_{\{a_n\}}$ es un $k_{\omega}$-grupo (ver  \cite[4.1.5]{ProZelLibro} ).
\end{itemize}
\end{remark}

\begin{Pro}\label{DualesCoincidenTopologicamente}

El grupo dual de $\Z_{\{p^n\}}$ coincide algebr\'aica y topol\'ogicamente con el dual de $(\Z,\lambda_p)$; es decir, se identifica con $\Z(\mathbf{p}^\infty)$ dotado de la topolog\'ia discreta.

\end{Pro}
\begin{Proof}
 En \cite[4.6]{PreprintDikran} se obtiene que el dual de $\Z_{\{p^n\}}$ coincide algebr\'aicamente con el dual de $(\Z,\lambda_p)$ y por tanto se identifica a $\Z(\mathbf{p}^\infty)$. La observaci\'on $(2)$ de \ref{kW} implica que $\Z_{\{p^n\}}^\wedge$ es metrizable y completo. Por el Teorema de Categor\'ia de Baire, teniendo en cuenta que  $\Z_{\{p^n\}}^\wedge$ es numerable, obtenemos que es discreto y por tanto coincide tambi\'en topol\'ogicamente con el dual de $(\Z,\lambda_p)$.

\end{Proof}

La proposici\'on \ref{DualesCoincidenTopologicamente} nos permite resolver negativamente un problema planteado en \cite{PreprintDikran}, que indicamos a continuaci\'on, aunque  para ello tenemos que aludir a una clase de grupos topol\'ogicos abelianos que en particular contiene a la clase de los grupos reflexivos.
 Se trata de  los grupos localmente cuasi-convexos, que constituyen una clase  de grupos que de alg\'un modo refleja y extiende las propiedades de los espacios vectoriales localmente convexos. De hecho todo espacio localmente convexo considerado como un grupo respecto de la adici\'on (es decir, olvidando la estructura lineal), es un grupo localmente cuasi-convexo. Estudiar a fondo esta clase de grupos implicar\'ia ampliar demasiado este trabajo, por tanto remitimos al lector interesado a \cite{tesisMontse} para formarse una idea clara y precisa de la misma. En \cite{PreprintDikran} se prueba - mediante c\'alculo directo de algunas envolturas cuasi-convexas- que los grupos  $\mathbb{Z}_{\{2^n\}}$ y $\mathbb{Z}_{\{3^n\}}$  no son localmente cuasi-convexos y se pregunta qu\'e ocurre para $\mathcal{T}_{\{p^n\}}$ con $p > 5$. De modo directo hemos obtenido lo siguiente:

\begin{Pro}

Sea $p$ un primo cualquiera. El grupo $\Z_{\{p^n\}}$ no es localmente cuasi-convexo y, por tanto, no es reflexivo.

\end{Pro}

\begin{Proof}
Para abreviar, llamemos $G : =\Z_{\{p^n\}}$. Si $G$ fuera localmente cuasi-convexo, la aplicaci\'on can\'onica $\alpha_G$ ser\'ia inyectiva y abierta (\cite[6.10]{tesislydia}). De nuevo por $(2)$ de la observaci\'on \ref{kW}, $G$ es k-grupo, lo que implica que  $\alpha_G$ es continua.  Por tanto  $\alpha_G$ es un encaje topol\'ogico de $G$ en $G^{\wedge\wedge}$. Por la proposici\'on \ref{DualesCoincidenTopologicamente}, $G^\wedge$ es discreto y en consecuencia  $G^{\wedge\wedge}$ es compacto. El mencionado encaje permite afirmar que  $G$ es precompacto. Ahora el  Teorema \ref{CR} nos indica que s\'olo hay una topolog\'ia precompacta en $\Z$ con dual  $\Z(\mathbf{p}^\infty)$, concretamente la p-\'adica $\lambda_p$. As\'i obtenemos la igualdad  $G=(\Z,\lambda_p)$, que contradice el  hecho de que $(\Z, \lambda_p)$ es metrizable pero  $G = \Z_{\{p^n\}}$ no lo es.

\end{Proof}
\begin{Quest}

\begin{itemize}
\item[(1)]?` Bajo qu\'e condiciones un par de $D$-sucesiones $\{a_n\}$ y $\{b_n\}$ dan lugar a grupos distintos $\Z_{\{a_n\}}$ y $\Z_{\{b_n\}}$?.
\item[(2)] Estudiar si la  modificaci\'on localmente cuasi-convexa de la topolog\'ia $\mathcal{T}_{\{a_n\}}$ es  de Mackey en el sentido definido en \cite{CMPT}.
 \end{itemize}

\end{Quest}

\end{document}